\documentclass[11pt]{article}
\usepackage{latexsym}
\usepackage{graphicx}
\usepackage{amsthm,amsmath}
\usepackage[OT2, T1]{fontenc}

\newcommand{\be}{\begin{equation} }

\newcommand{\ee}{\end{equation}}
\newcommand{\beq}{\begin{eqnarray} }
\newcommand{\eeq}{\end{eqnarray}}
\newcommand{\beqn}{\begin{eqnarray*} }
\newcommand{\eeqn}{\end{eqnarray*}}

\newcommand{\tends}{\rightarrow}

\def\dnav{\leavevmode\rlap{,\kern-0.075em,}\hphantom{''}}

\newcommand{\eps}{\varepsilon}

\newcommand {\R}{{\bf R}}
\newcommand{\Rp}{{\bf \overline{R}}}

\newcommand{\Med}{{\rm Med\, }} 
\newcommand{\E}{{\rm E\, }} 

\newcommand{\bol}[1]{{\mbox{\boldmath $#1$}}}
\newcommand{\cF}{{\cal F}}

\newcommand{\cC}{{\cal C}}
\newcommand{\cB}{{\cal B}}
\newcommand{\cA}{{\cal A}}

\newcommand{\cI}{{\cal I}}
\newcommand{\cJ}{{\cal J}}

\newcommand{\cH}{{\cal H}}
\newcommand{\cU}{{\cal U}}
\newcommand{\cV}{{\cal V}}
\newcommand{\ba}{{\bol{a}}}
\newcommand{\bb}{{\bol{b}}}
\newcommand{\bc}{{\bol{c}}}
\newcommand{\bd}{{\bol{d}}}
\newcommand{\bx}{{\bol{x}}}
\newcommand{\by}{{\bol{y}}}
\newcommand{\bX}{{\bol{X}}}
\newcommand{\bsX}{\mbox{\scriptsize\bf\em X}} 
\newcommand{\eop}{\hfill \raisebox{.3em}{\fbox{\rule{.2em}{0em}\rule{0em}{.2em}}}} 

\newcounter{secti}
\newcounter{item}[secti]

\numberwithin{equation}{section}

\theoremstyle{plain}
\newtheorem{thm}{Theorem}[section]
\newtheorem{lem}{Lemma}[section]
\newtheorem{cor}{Corollary}[section]

\theoremstyle{definition}
\newtheorem{exm}{Example}[section]
\newtheorem{dfn}{Definition}[section]

\begin{document}

\setcounter{page}{0} \thispagestyle{empty}

\vspace{5em}

\begin{center}

{\LARGE\bf Spatial medians, depth functions and multivariate
Jensen's inequality }


\bigskip

{\sc \DJ or\dj e Baljozovi\'c and Milan Merkle}

\bigskip

\parbox{25cc}{
{\bf Abstract. }{\small For any given partial order in a
$d$-dimensional euclidean space, under mild regularity assumptions,
we show that the intersection of closed (generalized) intervals
containing more than $1/2$ of the probability mass, is a non-empty
compact interval. This property is shared with common intervals on
real line, where the intersection is the median set of the
underlying probability distribution. So obtained multivariate
medians with respect to a partial order, can be observed as  special
cases of centers of distribution in the sense of type D depth
functions introduced by Y. Zuo and R. Serfling, {\em Ann. Stat.},
{\bf 28} (2000), 461-482. We show that the halfspace depth function
can be realized via compact convex sets, or, for example, closed
balls, in place of halfspaces, and discuss structural properties of
halfspace and related depth functions and their centers.  Among
other things, we prove that, in general, the maximal guaranteed
depth is $\frac{1}{d+1}$. As an application of these results, we
provide a Jensen's type inequality for functions of several
variables, with medians in place of expectations, which is an
extension of the previous work by M. Merkle, {\em Stat. Prob.
Letters}, {\bf 71} (2005), 277--281. }

}

\vspace{3em}

\end{center}

\vfill



\medskip


\noindent {\bf Acknowledgements.} This work is partially supported
by Ministry of Science and environmental protection of Serbia, under
the project number 144021. Second author acknowledges his partial
affiliation to Ra\v cunarski fakultet, Beograd, Serbia, and
Instituto de Mathematica, Universidade Federal do Rio de Janeiro,
Rio de Janeiro, Brasil.

\vfill

\noindent University of Belgrade, Faculty of Electrical Engineering,
P.O. Box 35-54, 11120 Belgrade, Serbia\\
emerkle$@$etf.bg.ac.yu\\
djordje.baljozovic$@$etf.bg.ac.yu

\newpage
\thispagestyle{plain}

\begin{center}

{\large\bf Spatial medians, depth functions and multivariate
Jensen's inequality}

\vspace{3em}

{\sc \DJ or\dj e Baljozovi\'c and Milan Merkle}

\end{center}

\medskip

\vspace{3em}

\section{Introduction.} To attain any median $m$ of a real valued
random variable $X$, we have to pass at least half of the
population, coming from either side of the real axis, via the
relations

\be \label{onemed}
 {\rm Prob} (X \in (-\infty, m])\geq p,\quad
   {\rm Prob} (X\in [m,+\infty)) \geq 1-p,\qquad p=\frac{1}{2}.
\ee

The proportion $\frac{1}{2}$ is the largest possible, in the
sense that there does not exist a $p>\frac{1}{2}$ such that
for arbitrary distribution of $X$ there exists $m\in \R$ which
satisfies (\ref{onemed}).  We  say that a median is the {\em deepest} point with respect to
a given distribution, or a data set. Quantitatively, we may assign a
{\em depth} to each point $x\in \R$, according to the value of a
{\em depth function} \be \label{deppth}
 D(x;P)= \inf \{ P((-\infty,x]), P([x,+\infty))\} ,\quad P(S):={\rm Prob} (X\in S),
 \ee
from where we can also see that the function $D$ attains its maximum
$1/2$ in the set of median points.

The literature devoted to the problem of extension of these
observations to higher dimensions contains a  variety of different
approaches - see \cite{small4} or \cite{zuoserf1} for a
comprehensive bibliography. In a high dimensional data set it is
desirable to select one point, or a set of points that would
correspond to intuitive notions of  "deepest point", "most central
point", or "the center of a data set", and can serve as a (best)
representative or a reference point of a distribution or a data set.

The definition (\ref{deppth}) of a depth function includes a notion
of direction: we can approach a point $x\in \R$ from either left or
right. The notion of direction in $\R^d$  is related to partial
order. In this paper we show, among other things, that  with respect
to any partial order that satisfies mild
 regularity conditions, there is a well defined median set, in the sense that
 the depth (to be defined in Section 3) of each point
of that median set is at least $1/2$. The median set is always
compact, and can be obtained as an intersection of a certain family
of sets.

The median sets defined with respect
to a given direction (i.e., a
partial order) are not affine invariant. On the other hand, Tukey's
median, or halfspace median (first introduced in \cite{tukey}, see also \cite{chentyl} or \cite{zuoserf1}),
which can be  defined as the set of maximal depth with respect to  the
depth function (commonly referred to as the halfspace depth)
 \be
\label{deptuk}
 D(x;P)= \inf \{P(H)\; |\; \mbox{$H$:\ any halfspace that contains $x$}\},
 \ee
  is affine invariant, but the maximal depth is not guaranteed to be $1/2$. As a generalization
 of (\ref{deptuk}), Zuo and Serfling in \cite{zuoserf1} offered a general "depth function of type D", which is
 defined as
 \be
 \label{deptyd}
 D(x;P, \cC)= \inf \{P(C)\; |\; x\in C\in \cC \},
 \ee
 where $\cC$ is a given collection of closed sets. A similar concept was introduced in \cite{small87}.
 In this paper we observe the function
 (\ref{deptyd}) defined with respect to arbitrary class $\cU$ of  open  sets
 in place of closed  sets in $\cC$. This slight change enables a considerable reduction of
 conditions for validity of certain results, and yields a collection of interesting and
 enlightening examples.   We prove several general
results for a class $\cU$ under  two very
 mild conditions, and we show that for each such class, the set of deepest points can be obtained as
 the intersection of a family of sets. It turns out that median sets can be also obtained as points
 having the depth not less than $1/2$, with a special choice of $\cU$.
 In this paper, by a multidimensional
 median set we understand  only the set of points with the depth not less than $1/2$. It can happen that
  such a set  is empty; in general case we use  use the term "center of the distribution" for the set of
  deepest points. These two cases need to be distinguished, as the corresponding sets have different properties.

  We show that in general, with the type D  depth
 in $\R^{d}$ we can be certain to find only   point(s)  with $D(x;P)\geq 1/(d+1)$, and we
 show that this bound is the best possible, thus extending  a result from \cite{dongas3} for
 the halfspace depth.

 Finally, we  show that each depth function generated by
a family  $\cU$ of complements of compact convex sets, can be
defined in terms of a family of halfspaces, and we represent the
halfspace center ("Tukey median") via intersection of level sets of
depth functions defined with respect to partial orders.
 In the last section, we discuss
 an analogue of Jensen's inequality for functions of several variables, with medians (or, generally,
 points in the center of distribution) in place of expectations.

The results presented in this paper are fully general, and hold
without any particular assumptions about the underlying
distribution.

For more details regarding Tukey's median see \cite{don1},
\cite{dongas3} and \cite{zuoserf1}. The latter paper and
\cite{liups1} give wide-ranging discussion of depth functions in
general with numerous examples and a summary of further researches
based on the halfspace depth.

In papers \cite{liu1,liu2} the simplical depth function was
proposed, primarily important for its corresponding version of
multivariate median; also, \cite{liu2} highlighted some features the
(simplical) depth function should fulfill.

Numerous other depth functions have been introduced in the
literature - some of the most relevant are presented in
\cite{nolan}, \cite{liusin1}, \cite{masstheo}, \cite{koshmos},
\cite{liups1}, \cite{zuoserf1,zuoserf4,zuoserf2}, \cite{varzhang},
\cite{zuo1}, \cite{zuo4}, \cite{zuocuiy}, \cite{zuo2} and references
therein. There are methods and techniques developed especially for
the data sets, see for example, \cite{chakchaud}. Robustness
properties of deepest points based on halfspace depth in context of
location statistics are studied in \cite{chentyl}.

Other approaches in defining multivariate centers are given in
\cite{small4} (many of them not related to depth functions),
\cite{niinoja}, \cite{varzhang}, \cite{hettrand} and references
therein; as can be seen, this concept is not at all unambiguous like
in the univariate case.

{\em Notations. } Throughout the paper, $\R$ denotes the set of real
numbers, and $\Rp$ is the set $\R$ together with $\pm \infty$. The
symbols $\R^d$ and $\Rp^d$ denote corresponding $d$-dimensional
euclidean spaces (Points in $\Rp^d$ are allowed to have $\pm \infty$
as coordinates). In Section 2, and in the material of later sections
related to intervals, the points in $\Rp^d$ with $d>1$ will be
denoted by bold letters $\bx, \by,\ldots $ and their coordinates by
$x_i, y_i, \ldots$, respectively. We say that the set $S\in \R^d$,
$d\geq 1$  is closed if it is closed in euclidean topology. Hence,
for $d=1$, the intervals $[a,+\infty)$ or $(-\infty, a]$ are closed
for any $a\in \R$; a similar remark holds for generalized intervals
with $d>1$. We will observe  $\R^d$-valued random variables $X$,
considering them as being measurable maps from some abstract
probability space $(\Omega, \cF, {\rm Prob})$ to $(\R^d, \cB^d, P)$.
Here $\cB^d$ is the Borel sigma-field on $\R^d$ and $P$ is a
probability measure induced on $\R^d$ by $X$ (the probability distribution
of $X$).

For a random variable
 $X$, $\Med X$ denotes any of its medians, and $\{ \Med X\}$ denotes the set of all medians. In the same way,
 we may talk about medians of a probability distribution $P$.

\medskip

\section{Multivariate medians.}
\label{mume}

We  start with a characteristic
property of univariate median set. Let $X$ be a random variable with
a probability distribution $P$ and let $J$ be any closed interval
with $P( J) >1/2$. We will show that $J$ contains every median of
$X$.  Indeed, if $m$ is a median of $X$ and
  $m\not\in J$,
then one of the intervals $(-\infty, m]$ or $[m,+\infty)$ is
disjoint with $J$, which is not possible, since the sum of
probabilities in both cases is greater than 1. Therefore, the
intersection
\[ \bigcap_{J=[a,b]:\ P(J) >1/2} J\]
is nonempty, and it contains the median set $[u,v]$ of $X$. Now,
observe that for  $J_{2n-1}=(-\infty, v+\frac{1}{(2n-1)}]$ and
$J_{2n} = [u-\frac{1}{2n}, + \infty)$, $n=1,2,\ldots$ we have that
$P(J_n)>1/2$ and so
\[ \{ \Med X \} = [u,v] = \bigcap_{n=1}^{+\infty}  J_n \supset \bigcap_{J=[a,b]:\ P(J) >1/2} J,\]
which together with the previous part, shows that \be \label{unimed}
 \{ \Med X \} =  \bigcap_{J=[a,b]:\ P(J) >1/2} J
 \ee

The relation (\ref{unimed}) can be as well taken as a definition of
the median set for a given distribution, and this definition can be
extended in a multidimensional environment if we choose one of many
possible extensions of the concept of one-dimensional interval. Out
of several ones that we may think of (convex sets, star-shaped sets,
balls and other special convex sets), only intervals with respect to
a partial order can do the work, to ascertain non-emptiness of the
intersection at the right hand side of (\ref{unimed}).

Let $\preceq$ be a partial order in $\Rp^d$ and let $\ba,\bb$ be
arbitrary points in $\Rp^d$. We define a $d$-dimensional
interval  $[\ba, \bb]$ as the set of points in $\R^d$ that are
 between $\ba$ and $\bb$:
\[ [ \ba, \bb] = \{ \bx\in \R^d\; |\; \ba \preceq \bx \preceq \bb\} \]

Note that the interval can be an empty set, or a singleton.
For the sake of simplicity, we want all intervals to
be topologically closed. The
interval can be norm bounded or norm unbounded; it would be
reasonable to expect intervals with finite "endpoints" to be norm
bounded, hence compact.
Further, we would expect that intervals can be "big" as we
wish, to contain any ball or any compact set. Finally, we expect
that bounded (with respect to partial order) sets posses the least
upper bound and greatest lower bound. To summarize, we assume the
following three technical conditions:
\begin{itemize}
\item[(I1)] Any interval $[\ba,\bb]$ is topologically closed, and
for any $\ba,\bb \in \R^d$ (i.e., with finite coordinates), the interval $[\ba,\bb]$ is a
compact set.
\item[(I2)] For any ball $B\subset \R^d$, there exist $\ba,\bb \in \R^d$
such that $B\subset [\ba, \bb]$.
\item[(I3)] For any set $S$ which is bounded from above with a finite point, there exists a finite
$\sup S$. For any set $S$ which is bounded from below with a finite
point, there exists a finite $\inf S$.
\end{itemize}

\begin{exm}
\label{ccorder} Let $K$ be a closed convex cone in $\R^d$, with
vertex at origin, and suppose that there exists a closed hyperplane
$\pi$, such that $\pi \cap K =0$ (that is, $K\setminus\{0\}$ is a
subset of one of open halfspaces determined by $\pi$). Define the
relation $\preceq$ by $\bx \preceq \by \iff \by -\bx \in K$. The
interval is then
\[ [\ba,\bb]= \{ \bx\; |\; \bx-\ba \in K \wedge \bb-\bx \in K\} = (\ba +K) \cap (\bb -K).\]
If the endpoints have some coordinates infinite, then the interval is
either $\ba +K$ (if $\bb \not\in \R^d$) or $\bb - K$ (if $\ba \not\in \R^d$) or $\R^d$ (if neither
endpoint is in $\R^d$).

It is not difficult to show that $\preceq$ is a partial directed
order, and that it satisfies conditions (I1)--(I3).

The simplest, coordinate-wise ordering, can be obtained with $K$ chosen to be the
orthant with $x_i \geq 0, i=1,\ldots , d$. Then

 \be
 \label{ocw}
 \bx \preceq \by \iff x_i \leq y_i ,\quad i=1,\ldots, d.
\ee

For the sake of illustration, let us note that possible kinds of intervals
with respect to the relation (\ref{ocw}) include:
\[ [(a_1,a_2), (b_1,b_2)],\ [(a_1,a_2), (b_1,+\infty)],\ [(a_1,a_2),(+\infty, b_2)],\]
\[ [(a_1,a_2),(+\infty,+\infty)],\ [(-\infty,-\infty), (b_1,b_2)],
\ [(a_1,-\infty), (+\infty, b_2)], \]
where  $a_1,a_2,b_1,b_2$ are real numbers. For infinite endpoints we use strict inequalities,
for example the last interval above is the set of $(x,y)\in \R^2$ such that $a_1\leq x <+\infty$ and
$-\infty <y \leq b_2$. The intervals may be empty; for example, the first listed interval is empty if
$a_1>b_1$ or if $a_2>b_2$.

\eop
\end{exm}

The following theorem extends the one-dimensional property  discussed in the beginning of
this section.

\begin{thm}
\label{jinters} Let $\preceq$ be a partial order in $\Rp^d$ such that
conditions (I1)--(I3) hold. Let $P$ be a probability measure on
$\R^d$  and let $\cJ$ be a family of intervals with respect to a
partial order $\preceq$, with the property that
\be \label{jprop}
P(J)
>\frac{1}{2},\quad\mbox{for each $J\in \cJ$}. \ee Then the
intersection of all intervals from $\cJ$ is a non-empty compact
interval.
\end{thm}

The compact interval claimed in the Theorem \ref{jinters} can be, in analogy to
(\ref{unimed}), taken as a definition of the median induced by the partial order $\preceq$:
 \be
\label{unimedd}
 \{ \Med \bX \}_{\preceq} \rule{1em}{0em} =  \bigcap_{J=[\ba,\bb]:\ P_{\mbox{\bsX}}(J) >1/2} J,
 \ee
 where $\bX$ is a random variable on $\R^d$ and $P_{\bsX}$ is its probability distribution.
 In what follows, we will omit the subscript if the underlying  relation $\preceq$ is obvious.

It is shown in Appendix (Lemma \ref{equalc}) that, in the case of the coordinate-wise
partial order,  the
characterization (\ref{unimedd}) is equivalent to a similar
characterization given in  \cite{repold}.
 It turns out that, in this
 case, the median set is just the Cartesian product of coordinate-wise median sets, which is
 the result stated in the form of an example in \cite{small87}.

 \begin{thm}
 \label{proj}
 Let the partial order $\preceq$ in $\R^d$, $d>1$,
 be defined by (\ref{ocw}), and let $\{\Med \bX\}$
 be the median set of a random vector $\bX\in \R^d$, defined with respect to the partial order $\preceq$.
 Then
 \be
 \label{cpcm}
 \{ \Med \bX \} = \{\Med X_1 \} \times \{ \Med X_2\}\times  \cdots \times \{ \Med X_d\}
 \ee
 \end{thm}

As we already mentioned, for other classes of sets, the intersection in (\ref{unimed}) is in general empty.
However, in the next section we will see that  some other classes  can also
have a non-empty intersection, but then, instead of having the measure $>1/2$, the sets will
generally have to have a greater measure.

\medskip

\section{Depth functions.} We start with a collection of  sets $\cV$ and the
collection $\cU$ that contains complements of sets from $\cV$, i.e.,
$ \cU = \{ S^c\; |\; S\in \cV\}$.
 For
each $x\in \R^d$ and any probability measure $P$ on Borel sets of
$\R^d$, define a depth function \be \label{dfg} D(x;P,\cU)= \inf \{
P(U)\; |\; x\in U\in \cU\}. \ee The role of $\cV$ will be clear
later in this section, when we  give an alternative description of
the depth function in terms of sets in $\cV$.

If for some $x$, there does not exist any $U\in \cU$ that contains
$x$, on the right hand side of (\ref{dfg}) we have empty set, and
then $D(x;P,\cU)=+\infty$. To avoid this, we assume that
\[ (C_1) \quad \mbox{  for every $x \in \R^d$ there is a $U\in \cU$ so that $x\in U$.}\]
Further, a constant depth function does not serve any purpose; to
avoid that situation, we may pose two additional conditions $(C_2)$
: \beqn &
(C_2')& \quad\mbox{
$D(x;P,\cU)>0$ for at least one $x\in \R^d$\qquad \mbox{and}}\\
& (C_2'') & \quad\lim_{\|x\|\tends +\infty}D(x;P,\cU) =0 \\
\eeqn
The condition $(C_2'')$ was also singled out in
\cite{zuoserf1}, as a requirement for any reasonable depth function.


Before proceeding further, let us see some examples.

\begin{exm}
\label{vared} $1^{\circ}$ The simplest family $\cU$ that satisfies
($C_1$) contains only one set - the whole space $\R^d$. Here
$D(x;P,\cU)=1$ for all $x$; condition $(C_2'')$ does not hold. In
next three examples, the conditions $(C_2')$ and $(C_2'')$ hold (see
Corollary \ref{cahd} in Appendix).

$2^{\circ}$ In $\R^d$, $d\geq 2$, let $\cV$ be a family of all
closed intervals $[\ba,\bb]$ with respect to the partial order
induced by a convex cone, as in Example \ref{ccorder}. It is easy to
see that for each point $x\in \R^d$, there exists a closed interval
$V$ such that $x\not \in V$; hence, the condition $(C_1)$ is
satisfied. A particular case of this example is the coordinate-wise
partial ordering in $\R^d$, which in $\R^2$ yields $V\in \cV$ to be
rectangles with sides parallel to axes.

$3^{\circ}$ Instead of intervals in previous examples, we can take
 $V\in \cV$ to be  arbitrary convex and compact sets
with a property that the collection $\cV$ is closed under
translations, and that sets $V$ can be arbitrary "big" (for example,
every ball in $\R^d$ should be contained in some $V\in \cV$).
Because of compactness and the translation property, for each $x\in
\R^d$ there exists a $V\in \cV$ that does not contain $x$, and
$(C_1)$ follows.

$4^{\circ}$ Consider now the class $\cV$ of all closed halfspaces.
The sets $U\in \cU$ are then open halfspaces, and then the
definition (\ref{dfg}) of depth function formally differs from
Tukey's halfspace depth function that requires closed halfspaces to
be in $\cU$. However, since $P(U)=\lim P(V_n)$, where $V_n\supset U$
are closed halfspaces with boundaries parallel to the boundary of
$U$ at euclidean distance $1/n$, it follows that values of $D$
coincide for these two cases. The condition $(C_1)$ is clearly
satisfied. \eop

\end{exm}

We are here interested chiefly in finding the set where the function
$D$ attains its global maximum, or, more generally, the sets of the
form

\be \label{salpha}
 S_{\alpha}= S_{\alpha}(P,\cU):=\{ x\in \R^d\; | \; D(x; P,\cU)\geq \alpha \} ,
\ee

The next Lemma gives a way to find $S_{\alpha}$ without  evaluation
of the depth function.

\begin{lem}
\label{depf1} Let $ \cU$ be any collection of non-empty sets in
$\R^d$, such that the condition $(C_1)$ holds. Then, for any
probability distribution $P$, \be \label{11a} S_{\alpha}(P,\cU) =
\bigcap_{V\in \cV, P(V)>1-\alpha} V, \ee for any $\alpha \in (0,1]$
such that there exists a set $U\in \cU$ with $P(U)<\alpha$;
otherwise $S_{\alpha}=\R^d$.
\end{lem}

It is instructive  first to observe $S_{\alpha}$ in $d=1$, as in the
next example.

\begin{exm}
Let $X$ be a real random variable with the distribution $P$. Take
$\cV$ to be the family of all  closed intervals in $\R$,  and $\cU$
to be the family of their complements. Then, using similar arguments
as in the proof of (\ref{unimed}) at the beginning of Section 2, it
can be derived that $S_{\alpha}=[q_{\alpha},Q_{1-\alpha}]$, where
$q_{\alpha}$ is the smallest quantile of $X$ of order $\alpha$, and
$Q_{1-\alpha}$ is the largest quantile of $X$ of order $1-\alpha$:

\beq \label{qab}
 q_{\alpha} &=&\inf\{ t\in \R \; |\; {\rm Prob}\, ( X\leq t) \geq \alpha\}\quad\mbox{ and}\nonumber \\
Q_{1-\alpha} &=&\sup\{ t\in \R \; |\; {\rm Prob}\, (X \geq t) \geq
\alpha\}. \eeq

For $\alpha=\frac{1}{2}$,  $[q_{\frac{1}{2}},Q_{\frac{1}{2}}]$ is
the median interval.\eop

\end{exm}

In the case when  the family $\cV$ is consisted of closed intervals
with respect to a partial order that satisfies (I1)--(I3),  it
follows from  Section 2 and Lemma \ref{depf1}, that the set
$S_{1/2}$ is non-empty, i.e. the corresponding depth function has
maximum which is $\geq 1/2$, regardless of distribution $P$. For
other families of $\cV$, the guaranteed maximum is smaller.

\begin{exm}
\label{e13} Consider the halfspace depth, as in Example $4^{\circ}$
of \ref{vared}, in $\R^2$, with the probability measure $P$ which
assigns mass $1/3$ to points $A(0,1)$, $B(-1,0)$ and $C(1,0)$ in the
plane. Each point $x$ in the closed triangle $ABC$ has
$D(x)=\frac{1}{3}$; points outside of the triangle have $D(x)=0$.
So,  the function $D$ reaches its maximum value $\frac{1}{3}$.

Let us now observe the same distribution, but with depth function
defined with the family $\cV$ of closed disks. The intersection of
{\em all} closed disks $V$ with $P(V)>2/3$ is, in fact, the
intersection of all disks that contain all three points $A,B,C$, and
that is the closed triangle $ABC$. For any $\eps>0$, a disc $V$ with
$P(V)> 2/3-\eps$ may contain only two of points $A,B,C$, but then it
is easy to see that the family of all such discs has the empty
intersection. Therefore, $S_{\alpha}$ is non-empty for $\alpha\leq
1/3$, and again, the function $D$ attains its maximum value $1/3$ at
the points of closed triangle $ABC$.

If the depth function is defined in terms of rectangles with sides
parallel to coordinate axes, then the maximum depth is $2/3$ and it
is attained at $(0,0)$. This conclusion follows immediately from
Theorem \ref{proj}.\eop

\end{exm}

If $\alpha_m$ is the maximum value of $D(x;P,\cU)$ for a given
distribution $P$, the set $S_{\alpha_m}$, i.e., the set of deepest
points with respect to $P$, is called the center of the distribution
$P$, and will be denoted by $C(P,\cU)$.

In the next theorem, we discuss some properties of the center, in
the case when sets in $\cU$ are open. A similar  result for the
family $\cU$ of closed sets was obtained in \cite[Theorem 2.11]{zuoserf1},
but under more restrictive assumptions.


\begin{thm}
\label{maint} Let $\cV$ be a collection of closed subsets of $\R^d$,
and let $\cU$ be the collection of  sets $V^c$, $V\in
\cV$,  such that the condition $(C_1)$ holds. Then, for arbitrary
probability measure $P$, the function $x\mapsto D(x; P, \cU)$ is
upper semicontinuos. Moreover, under conditions $(C_2)$, the set
$C(P,\cU)$ on which $D$ reaches its maximum is equal to the minimal
nonempty set $S_{\alpha}$, that is,
\[ C(P,\cU)=\bigcap_{\alpha:S_\alpha \neq \emptyset} S_{\alpha} (P,\cU).\]
The set $C(P,\cU)$ is a non-empty compact set and it has the following
representation: \be \label{cedis} C(P,\cU)= \bigcap_{V\in \cV,
P(V)>1-\alpha_m} V ,\quad \mbox{where $\alpha_m= \max_{x\in \R^d}
D(x;P,\cU)$.} \ee
\end{thm}

\medskip

\section{Equivalence of depth functions.} We already noticed that the
depth functions with $\cU$ being all open or all closed halfspaces,
have the same value at every point. So, it is possible that two
different classes of sets in place of $\cU$ generate the same depth
function.  In the next theorem we give a sufficient condition for
the equivalence of two depth functions.

\begin{thm}
\label{equco} Let $\cA$ and $\cB$ be families of subsets of $\R^d$.
Suppose that the condition $(C_1)$ holds for at least one of these
families, and, in addition, the following condition $(E)$:

\begin{itemize}
\item[$(E')$] For each $A\in \cA$, $A=\bigcup_{B\in \cB, B\subset A} B$, and
\item[$(E'')$] For each $B\in \cB$, there exists at most
 countable collection of sets $A_i \in \cA$, such that
$A_1\supseteq A_2\supseteq \ldots$ and $B=\bigcap_{i} A_i$.
\end{itemize}

 Then the condition $(C_1)$ holds for both $\cA$ and $\cB$ and
 depth function with respect to both families  are equal, with any probability distribution $P$:
 \[ D(x;P,\cA)=\inf\{ P(A)\; |\; x\in A\in \cA \} = \inf\{ P(B)\; |\; x\in B\in \cB \} = D(x; P, \cB)\]
\end{thm}

 An important application of Theorem \ref{equco} is to establish the equivalence
of depth functions defined by a family of open sets $A\in \cA$ and
their topological closures $\bar{A}\in \cB$.
 In this setup, we note that $(E)$  holds  in cases when $\cA$ is invariant with respect to
 translations, and
consists of  (i) open halfspaces or (ii) complements of  all closed
intervals with any convex cone partial order,
 as in example \ref{ccorder}. Another application of Theorem \ref{equco} will be
 given in Theorem \ref{tequch}.

\medskip

\section{Convex sets and halfspaces.}

Suppose that a family $\cU$ contains a sequence of nested sets $U_n$ that intersect
at one point $x\in \R^d$. Then  $D(x; P,\cU)=P(x)$ for any $P$,  which is undesirable
property. Therefore $\cU$ should not contain sets that shrink to a
point. A way to avoid that is to choose  sets in $\cU$ to be
unbounded, or to choose sets in $\cV$ to be bounded.

Further, it is natural to have a convex center of distribution,
which is achieved (via Theorem \ref{maint}) if sets in $\cV$ are
convex. With  non-convex sets in $\cV$, and with a discrete
distribution $P$, we can again have that $D(x;P,\cU)=P(x)$ for every
$x\in \R^d$, as shown in the next example.

\begin{exm}
\label{e1g}

In $\R^2$, let $K$ be the lower half of the first quadrant,  bounded
by halflines $y=0$ and $y=x$. Let us consider the family ${\cU} $ of
sets that can be obtained by  arbitrary rotations and
translations of $K$.  Let $\cV$ be the family of complements of sets in
$\cU$: sets in $V\in \cV$ are non-convex.


Suppose that a distribution $P$ is concentrated in six points
$X_{1,2}(\pm 1,0)$, $X_{3,4}(\pm 2, 0)$, $X_{5,6}(0,\pm 1)$ (as in
the Counterexample 2 in \cite{zuoserf1}). Then for each point $x\in
\R^d$, there is a set $U\in \cU$ that does not contain any point
of the support different from $x$;
hence the depth of each point is $P(\{x\})$. With a specific discrete
distribution $P$, the center is the point with the greatest
probability mass. If $P$ is uniform across the $X_i$, $i=1,\ldots,
6$, then the center is the discrete set $\{ X_1,\ldots, X_6\}$. This
is in a sharp contrast with the halfspace depth function, which in
this case yields the single point center at $(0,0)$, with the
maximum depth $\frac{1}{2}$.\eop

\end{exm}

A prototype of  depth functions that we discuss in this section is a
depth function defined  with respect to families $\cU$ of
complements of compact convex sets. In the light of the arguments
given above, these requirements are natural and they are not too
restrictive (see also Lemma \ref{condit}). Although it may look that
by these requirements we are excluding the halfspace depth from
consideration, it is not so, as we will see after the Theorem
\ref{tequch}.

From the material of Section 2, it follows that the  depth function
based on a family $\cV$ of intervals, attains the maximal value of
at least $1/2$, regardless of the dimension $d$. In general, the
maximum depth with a family $\cV$ of convex sets, can not be smaller
than $\frac{1}{d+1}$. This conclusion follows from the next theorem,
which is an extension of results in \cite{dongas3} and \cite{RouRut}.

\begin{thm}
\label{lept} Let $P$ be any probability measure on Borel sets of
$\R^d$. Let $\cV$ be any family of closed convex sets in $\R^d$, and
let $\cU$ be the family of their complements. Assume that conditions
$(C_1)$ and $(C_2'')$ hold. Then the condition $(C_2')$ also holds,
and there exists a point $x\in \R^d$ with $D(x; P,\cU) \geq
\frac{1}{d+1}$.
\end{thm}

The lower bound for $D$ in Theorem \ref{lept} is the greatest
generally possible. As the next example shows, for the halfspace
depth, in any dimension $d\geq 1$, there exist a probability measure
$P$ such that $D(x; P,\cU) \leq \frac{1}{d+1}$ for all $x\in \R^d$.

\begin{exm}
This is an extension of the example \ref{e13}. Let $A_1,\ldots,
A_{d+1}$ be points in $\R^d$ such that they do not belong to the
same hyperplane (i.e. to any affine subspace of dimension less than $d$),
and suppose that $P(\{ A_i\}) = \frac{1}{d+1}$ for each
$i=1,2,\ldots, d+1$. Let $S$ be a closed $d$-dimensional simplex
with vertices at $A_1,\ldots, A_{d+1}$, and let $x\in S$. If $x$ is
a vertex of $S$, then there exists a closed halfspace $H$ such that
$x\in H$ and other vertices do not belong to $H$; then
$D(x)=P(H)=1/(d+1)$. Otherwise, let $S_x$ be a $d$-dimensional
simplex with vertices in $x$ and $d$ points among $A_1,\ldots,
A_{d+1}$ that make together an affinely independent set. Then for
$S_x$ and the remaining vertex, say $A_1$,  there exists a
separating hyperplane $\pi$  such that $\pi \cap S_x = \{x\}$ and
$A_1\not\in \pi$ (see \cite[Section 11]{rocka}). Let $H$ be a
halfspace with boundary $\pi$, that contains $A_1$. Then also
$D(x)=P(H)=1/(d+1)$.  So, all points  $x\in S$ have $D(x)=1/(d+1)$.
Points $x$ outside of $S$ have $D(x)=0$, which is easy to see. So,
the maximal depth in this example is exactly $1/(d+1)$. \eop

\end{exm}

In fact, if we have a family of compact convex sets $\cV$ that
contain  arbitrary large sets (in the sense of the following lemma),
then it is sufficient to assume only condition $(C_1)$, and then
$(C_2)$ will automatically hold. A natural way to choose $\cV$ would
be then, to choose one compact convex shape, and allow translations
(and, possibly, rotations, if we want an affine invariant depth).

\begin{lem}
\label{condit} Let $\cV$ be a family of compact convex sets in
$\R^d$, and let $\cU$ be the family of complements of sets in $\cV$,
such  that the condition $(C_1)$ holds. Suppose that for every
closed ball $B\in \R^d$ there exist a  set $ V \in \cV $, such that
$B\subset V$.  Then the family $\cU$ and the depth function
$D(\cdot\ ;P,\cU)$ satisfy conditions $ (C_2')$ and $(C_2'')$, with
any probability measure $P$ on $\R^d$.
\end{lem}

In the next theorem, we use the fact that every closed convex set
can be represented as an intersection of closed halfspaces (see, for
example, \cite[Theorem 11.5]{rocka}). This representation is not unique (and we do not need
uniqueness neither in the statement nor in the proof); however,
there is a unique minimal representation of a convex set as the intersection
of all its tangent halfspaces \cite[Theorem 18.8]{rocka}, which is an intuitive
model for the representation (\ref{tequch1}) below.

\begin{thm}
\label{tequch} Let $\cV$ be a collection of closed convex sets and
$\cU$ the collection of complements of all sets in $\cV$. For each $V\in
\cV$, consider a representation
\be
\label{tequch1}
 V = \bigcap_{\alpha\in A_{V}}   H_\alpha,
 \ee
where $H_{\alpha}$ are closed subspaces and $A_V$ is an index set. Let
\[ \cH^V = \{ \overline{H^c_{\alpha}}+x\; |\; \alpha \in A_V,\ x\in \R^d\} \]
be the collection of closures of complements of halfspaces $H_{\alpha}$ and their translations.
Further, let
\[ \cH
=\bigcup_{V\in \cV} \cH^V. \]
If for any $H\in \cH$ there exists at most countable collection of sets
$V_i\in \cV$, such that \be \label{tequch2}
 V_1\subseteq V_2\subseteq \cdots\quad  \mbox{and}\quad  \stackrel{\circ}{H}= \bigcup V_i,
 \ee
then
\[ D(x; P,\cU)=D(x; P, \cH)=D(x; P, \stackrel{\circ}{\cH}),\quad \mbox{ for every $x\in \R^d$},\]
where $\stackrel{\circ}{\cH}$ is the family of open halfspaces from
$\cH$.
\end{thm}

As a corollary to Theorem \ref{tequch}, we can single out two important particular cases.
Conditions (\ref{tequch1}) and (\ref{tequch2}) in both cases can be easily proved.

\begin{cor}
\label{cortequch}
{\bf a)} Let $\cV$ be the family of closed intervals with respect to the
partial order defined with a convex cone $K$, as in the Section \ref{mume}. Then
for any probability distribution and any $x\in \R^d$,
\[ D(x; P,\cU)=D(x; P, \cH),\]
where $\cU$ is the family of complements of sets in $\cV$ and $H$ is the family of all tangent halfspaces
to $K$, and their translations.

In particular, if $\cV$ is the family of intervals with respect to the  coordinate-wise partial order, then
the corresponding depth function is the same as the depth function generated by halfspaces with borders parallel
to the coordinate hyperplanes.

{\bf b)} Let $\cH$ be the family of all closed halfspaces, let $\cU_k$ be the family of complements of all
compact closed sets, and let $\cU_b$ be the family of complements of closed balls in $\R^d$. Then

\[ D(x; P, \cH)=D(x; P,\cU_k)=D(x; P,\cU_b),\]

 That is, the Tukey halfspace depth can be realized via complements of closed convex sets or via complements of
closed balls.

\end{cor}

The second part of Corollary \ref{cortequch} implies, via Lemma \ref{depf1}, that for the
halfspace depth function $D$, we have \be \label{salphad}
S_{\alpha}=\{ x\in \R^d\; | \; D(x)\geq \alpha \} = \bigcap_{K:\
P(K)>1-\alpha} K = \bigcap_{B:\ P(B)>1-\alpha} B, \ee where $K$ are
compact convex sets, and $B$ are closed balls. The reduction to balls is of the obvious interest in
applications, where we have to find deepest points of a high dimensional cloud of data.

\section{Affine invariance and another representation of the
halfspace center of distribution.} A depth function $D(x;P_X,\cU)$ in
$\R^d$ is said to be affine invariant, if
\be \label{ain}
 D(Ax +b; P_{AX+b}, \cU)=D(x;P_X,\cU)\quad\mbox{for any probability measure $P$},
 \ee
for any nonsingular $d\times d$ matrix $A$, any $b\in \R^d$ and
$x\in \R^d$, where $P_{AX+b}$ is a probability distribution of a
random variable $AX+b$, $X$ being a random variable with the
distribution $P_X$. From the definition (\ref{ain}), it follows that
one sufficient condition for affine invariance of $D$ is the affine
invariance of $\cU$: If $U\in \cU$, then $AU+b\in \cU$, for all $A$
and $b$. This condition is satisfied with the family $\cU$ of all
halfspaces; hence the halfspace depth is affine invariant. Due to
the fact that the same depth function can be generated by different
families $\cU$, this condition is not necessary, as the next example
shows.

\begin{exm}
Let $\cV$ be the family of all closed discs in $\R^2$, and let $\cU$
be the family of their complements. The family $\cU$  is not affine
invariant, because the circles
 transform into ellipses, with a non-orthogonal
matrix $A$. However, the family of all halfplanes $\cH$ generates
the same depth function as the family $\cU$, and so, the depth $D(x;
P,\cU)$ is equivalent to halfspace depth, hence, it is affine
invariant.\eop

\end{exm}

For depth functions that can be generated by a family of halfspaces,
the conditions of affine invariance  can be expressed via
translation and rotations, as every halfspace in $\R^d$ can be
transformed into another one by one rotation and one translation.
That is, for every two halfspaces $H_1,H_2\in \R^d$, there exists an
affine transformation $x\mapsto Ax+b$ with $A$ being an orthogonal
matrix, such that $H_2= AH_1+b$.

Consider one coordinate system in $\R^d$, with the  corresponding
set $\cJ$ of coordinate-wise intervals.
 Any rotation $\rho$
of the coordinate system will produce another family of intervals
$\cJ_{\rho}$. According to Theorem \ref{tequch}, depth functions
based on the family $\cup_{\rho}\cJ_{\rho}$ (where the union goes
through all possible rotations) is equivalent to the halfspace depth
function. More generally, we may observe any set of partial orders
$\{\preceq_{\rho}\}$ (where $\rho$ belongs to some index set) that
satisfy conditions (I1)-(I3) of Section 2 such that the
corresponding families of intervals $\cJ_{\rho}$ (i.e., families
$\cU_\rho$ of complements of sets from $\cJ_\rho$) together generate
the halfspace depth function. Let us call  such set of partial
orders  {\em complete}. For a given probability distribution $P$,
and a complete set of partial orders, let
\[ S_{\alpha,\rho} = \{ x\in \R^d\; | \; D(x; P, \cU_{\rho})\geq \alpha \} =
 \bigcap_{J\in \cJ_{\rho}:\ P(J)>1-\alpha} J;\]
\[ S_{\alpha}=\{ x\in \R^d\; | \; D(x; P, \cH)\geq \alpha \} ,\]
where $\alpha \leq \alpha_{m}$, and $\cH$ is the family of all open
halfspaces. Then by completeness, we have that \be \label{cendis}
 S_{\alpha}=\bigcap_\rho \bigcap_{J\in \cJ_{\rho}:\ P(J)>1-\alpha} J = \bigcap_{\rho} S_{\alpha,\rho}.
 \ee
For $\alpha=\alpha_m$, (\ref{cendis}) gives another representation
of the center of a distribution, in terms of sets that are not
affine invariant.  If we take any finite subset of partial orders,
$\rho=1,\ldots, n$, then we have \be \label{css}
 C(P,\cH) \subset \bigcap_{\rho=1}^n S_{\alpha_m,\rho},
 \ee
which gives an upper bound for the center of distribution in terms
of finitely many partial orders.

Note that, in general, the sets $S_{\alpha_m,\rho}$ are not centers
of the distribution with respect to $\preceq_{\rho}$; we proved in
Section 2 that there exist median sets $S_{1/2,\rho}$. In general,
we have that $\alpha_m <1/2$, and $S_{1/2,\rho}\subset
S_{\alpha_m,\rho}$. Median sets with respect to different partial
orders may have empty intersection. For example, if the distribution
is absolutely continuous, then every median set with respect to a
coordinate-wise partial order is a singleton; clearly by a rotation
of the coordinate system we may obtain different singletons.

%
%

\medskip

\section{A version of Jensen's inequality.} Let $\cV$ be a family of
closed sets, $\cU$ the family of complements of sets from $\cV$ and
let $D(x; P, \cU)$ be defined as in previous sections. Let
$C=C(P,\cU)$ be the center of a probability measure
 $P$
in $\R^d$, with $\alpha_m$ being the maximum of the depth function.
Assume conditions $(C_1)$ and $(C_2)$.

For a random variable $X$ with the distribution  $P$, the points in
the set $C(P,\cU)$ can be thought of as a kind of mean values of
$X$, in the same sense as univariate medians are being thought of.
If $f$ is a real valued function defined on $\R^d$, then the
analogous mean value of $f(X)$ are points in the closed interval
$[q_{\alpha_m},Q_{1-\alpha_m}]$, where, by (\ref{qab}),
$q_{\alpha_m}$ is the smallest quantile of $f(X)$ of order
$\alpha_m$, and $Q_{1-\alpha_m}$ is the largest quantile of order
$1-\alpha_m$. If $\alpha_m = 1/2$,   we have the  median interval of
$f(X)$, and the center $C(P,\cU)$ becomes $\{ \Med X\}$.  Let $m\in
\{ \Med X\}$ and $M\in \{\Med f(X)\}$. With analogy to Jensen's
inequality $f(\E X)\leq \E f(X)$ for convex functions, we may expect
that $f(m)\leq M$ for an appropriate class of functions $f$. Indeed,
we prove a  result of that kind,   for the class of functions that
are described in the following definition. The name {\em C-function}
is taken from \cite{repold}, where it was used in a more particular
context.

\begin{dfn}
\label{cfun} A function $f: \R^d\mapsto \R$ will be called a
\emph{C-function} with respect
 to a given  family $\cV$ of closed subsets of $\R^d$, if
 $ f^{-1}((-\infty,t])\in \cV$ or is empty set, for every $t\in \R$.
\end{dfn}

 \begin{exm}
 \label{excf}

 $1^{\circ}$ If $\cV$ is the family of all closed convex sets in $\R^d$,
 then the class of
  corresponding  C-functions is precisely the class of lower continuous
  quasi-convex functions, i.e., functions $f$
  that have the property that $ f^{-1}((-\infty,t])$ is a closed set for any $t\in \R$ and
  \[
   f(\lambda x +(1-\lambda)y) \leq \max \{ f(x), f(y)\} ,\qquad \lambda\in [0,1],\quad x,y\in \R^d.
   \]
  This is easy to see, starting from the definition \ref{cfun}. In particular, every convex
  function on $\R^d$ is a C-function with respect to the class of all convex
  sets.

 $2^{\circ}$  Let $D(x)$ be a halfspace depth function. Then it follows from (\ref{11a}) that the sets
  $S_{\alpha}$ are convex, which implies that the function $x\mapsto 1-D(x)$ is a C-function
  with respect to a family of all closed convex sets.

$3^{\circ}$A function $f$ is a  C-functions with respect to a family
of closed intervals (with respect to some partial  order), if and
only if
\[ \{ \bx\in \R^d \; |\; f(\bx) \leq t \} = [\ba,\bb ] ,\qquad \mbox{for some $\ba,\bb \in \R^d$ }, d\geq 1.\]
This condition is not satisfied for all convex functions. For
example, in $\R^2$, with the coordinate-wise partial order the
function defined by $f(x,y)= x^2 +y^2$ is convex, but the sets $\{
(x,y)\; |\; x^2+y^2\leq t\}$ are not intervals.

$4^{\circ}$ In $\R^2$, with coordinate-wise intervals, the function
$f$ defined by
\[ f(x,y)= \max \{ |x-a_1|-|x-b_1|,\ |y-a_2|-|y-b_2|\} \]
is a C-function, where $\ba (a_1,a_2)$ and $\bb (b_1,b_2)$ are given
points in $\R^2$.

$5^{\circ}$ In general, assuming conditions (I1)-(I3), we may define
the depth function $D(x)$ with respect to the class $\cU$ of
complements  of the given family of intervals. Since the
intersection of closed intervals is again a closed interval, we see
that here also the function $x\mapsto 1-D(x)$ is a C-function.

$6^{\circ}$ Note that, since we require sets in $\cV$ to be closed,
every C-function is lower semicontinuous.\eop

\end{exm}

The next two theorems are versions of Jensen's inequality.

\begin{thm}
\label{jenmedia} Let $\cV$ be a family of closed intervals with
respect to a partial order in $\R^d$, such that conditions
(I1)--(I3) are satisfied. Let $\{ \Med X\}$ be the median set of a
random variable $X$ with respect to the chosen partial order, and
let $f$ be a C-function with respect to the family $\cV$. Then for
every $M\in \Med \{f(X)\}$, there exists an $m\in \{\Med X\}$, such
that \be \label{jenmi}
 f(m) \leq M.
\ee In particular, if $m$ or $M$ are unique, then (\ref{jenmi})
holds for any $m, M$.
\end{thm}

In general case the depth function does not necessarily reach the value of
$1/2$, and we have only a weaker result:

\begin{thm}
\label{jenmedib} Let $\cV$ be a family of closed subsets of $\R^d$,
and let $\cU$ be the family of their complements. Assume that
conditions $(C_1)$ and $(C_2)$ hold with a given probability measure
$P$, induced by a random variable $X$. Let $\alpha_m$ be the maximum
of the depth function $D(x;P,\cU)$,  which is achieved in all points
of the center $C(P,\cU)$ and let $f$ be a C-function with respect to
$\cV$. Then for every $m\in C(P,\cU)$ we have that \be
\label{jenmig}
 f(m)\leq Q_{1-\alpha_m},
\ee where $Q_{1-\alpha_m}$ is the largest quantile of order
$1-\alpha_m$ for $f(X)$.
\end{thm}

To show that we can not claim anything better in a general case,
consider the following example:

\begin{exm}
Let $A, B, C$ be non-colinear points in the two dimensional plane,
and let $\cH$ be the collection of closed halfplanes. Let
$l(AB),l(AC),l(BC)$ be the lines determined by two indicated points.
Let $H_1$ be the closed halfspace that does not contain the interior
of the triangle $ABC$ and has  $l(AB)$ for its boundary, and let
$H_2$ be its complement. Define a function $f$ by
\[ f(x)= e^{-d(x,l(AB))} \quad\mbox{if $x\in H_1$ },\qquad f(x)= e^{d(x,l(AB))}\quad\mbox{if $x\in H_2$},\]
where $d(\cdot,\cdot)$ is euclidean distance.
 Then $f(A)=1$, $f(B)=1$ and $f(C)>1$, and $f$ is clearly a
$C$-function with respect to the class $\cH$. Now suppose that $P$
assigns mass $1/3$ to each of the points $A,B,C$. Then, by example
\ref{e13}, we know that the center $C(p,\cH)$ of this distribution
is the set of points of the triangle $ABC$, with $\alpha_m=1/3$.
 Hence, for $m\in C(P,\cH)$, $f(m)$ takes all values in $[1, f(C)]$.
 On the other hand, quantiles of the order $2/3$ are points in
 the closed interval $[1, f(C)]$; hence the most we can state is that $f(m)\leq f(C)$, with $f(C)$ being
 the largest quantile of order $2/3$.

\end{exm}

\medskip

\appendix

\begin{center}
\section{APPENDIX: PROOFS AND AUXILIARY RESULTS}\label{app}
\end{center}

\refstepcounter{secti}

In order to prove Theorem \ref{jinters}, we need the following
lemma.

\begin{lem}
\label{lemma1} Let $\preceq$ be a partial order in $\R^d$ such that
the conditions (I1) and (I3) hold. Let
\[ \cJ=\{ J_{\alpha}\; |\; J_{\alpha}=[\ba^{\alpha} ,\bb^{\alpha} ],\quad \alpha \in A\}\]
 be a collection of closed intervals, where $A$ is an index set. Assume that there is at least one
 $\alpha$ such that $\ba^{\alpha}\in \R^d$ (i.e., have all coordinates
 finite) and at least one $\beta$ such that $\bb^{\beta}\in \R^d$. Suppose that $J_{\alpha} \cap
 J_{\beta} \neq \emptyset$ for all $\alpha, \beta$. Then
 \begin{itemize}
 \item[(i)] $\ba^{\alpha} \preceq \bb^{\beta}$, for any $\alpha,\beta \in A$;
 \item[(ii)] The intersection of all sets in $\cJ$ is a non-empty compact interval $[\ba,\bb]$, with
 $\ba,\bb\in \R^d$.
 \end{itemize}
\end{lem}

\begin{proof}
If intervals $[\ba, \bb]$ and $[\bc, \bd]$ have a common point
$\bx$, then $\ba \preceq \bx \preceq \bb$ and $\bc \preceq\bx
\preceq \bd$; hence $\ba \preceq \bd$ and $\bc\preceq \bb$. This
shows (i). Further, to show (ii), note that by assumptions and (i),
the set $\{ \ba^{\alpha}, \ \alpha \in A\}$ is bounded from above
with a finite point, and so by (I3), there exists
$\ba=\sup_{\alpha\in A} \ba^{\alpha}$. In an analogous way we
conclude that there exists $\bb=\inf_{\beta\in A} \bb^{\beta}$. By
properties of the infimum and supremum, we have that
$\ba^{\alpha}\preceq\ba\preceq \bb\preceq\bb^{\alpha}$, for all
$\alpha\in A$, so the interval $[\ba,\bb]$ is non-empty; it is
compact by assumption (I1), and it is contained in all intervals of
the family $\cJ$. On the other hand, any point $\bc$ that is common
for all intervals $J_{\alpha}$ must be an upper bound for $\{
\ba^{\alpha}\}$ and a lower bound for $\bb^{\alpha}$; hence
$\ba\preceq \bc \preceq\bb$, that is, $\bc \in [\ba,\bb]$, and (ii)
is proved.
\end{proof}

{\sc Proof of the Theorem \ref{jinters}.} It is clear that any two
intervals in $\cJ$ have a non-empty intersection; besides, by (I2),
at least one of the intervals has finite endpoints. Then the
assertion follows by Lemma \ref{lemma1}. \eop

{\sc Proof of the Theorem \ref{proj}.} To simplify notations, we
give the proof for $d=2$; the proof in a general case is analogous.
We have the sequence of relations \beqn && \{ \Med \bX \} =
 \bigcap_{[\ba,\bb]:\ {\rm Prob} (\bX \in [\ba,\bb])>\frac{1}{2}} [\ba,\bb]\\
      &=&  \bigcap_{[\ba,\bb]:\ {\rm Prob} (\bX \in [\ba,\bb])>\frac{1}{2}}[a_1,b_1]\rule{2em}{0em}\times
            \bigcap_{[\ba,\bb]:\ {\rm Prob} (\bX \in [\ba,\bb])>\frac{1}{2}}[a_2,b_2]\\
      &\supset& \bigcap_{[a_1,b_1]:\ {\rm Prob} (X_1 \in [a_1,b_1])>\frac{1}{2}} [a_1,b_1]\rule{2em}{0em}\times
            \bigcap_{[a_2,b_2]:\ {\rm Prob} (X_2 \in [a_2,b_2])>\frac{1}{2}}   [a_2,b_2])\\
            &=& \{ \Med X_1\} \times \{  \Med X_2\},
\eeqn where we used the fact that $X_1\in [a_1,b_1]$ whenever $\bX
\in [\ba,\bb]$. On the other hand, if $\{ \Med X_1\} =[a,b]$ and $\{
\Med X_2\}=[c,d]$, then we have that
\[ [a,b]\times [c,d] = [a,+\infty)\times \R \cap
(-\infty,b] \times\R \cap \R\times  [c,+\infty)\cap \R\times
(-\infty, d]\] and we note that all four two-dimensional intervals
on  the right hand side of the last identity, can be expressed as
intersections of a sequence of intervals $J_n$ with {\rm Prob}$(\bX
\in J_n)>1/2$; for example,
\[ [a,+\infty)\times \R = \bigcap_{i=n}^{+\infty} [a-1/n,+\infty)\times \R ,\]
and ${\rm Prob} (\bX \in [a-1/n,+\infty)\times \R) >1/2$ because
$[a,b]$ is the median set for $X_1$. From this we conclude that
\[ \{ \Med X_1 \}\times \{ \Med X_2\} \supset \{ \Med \bX\} ,\]
and the theorem is proved. \eop

Another result related to the coordinate-wise partial order intervals is
presented in the next lemma. For a given probability measure $P$,
denote by $\cJ$ the class of all intervals $J$ with the
property that $P (J)>1/2$.

In \cite{repold}, the class
$\cI$ is defined as the family of all  closed intervals $I\subset \R^d$ (with respect
to coordinate-wise partial order) with the following property: If $J$ is any closed
interval (with respect to the same partial order) that contains
$I$ as a proper subset, then $J\in \cJ$.

In the following lemma, we  show that the intersection
of all intervals in $\cI$ coincides with the intersection of all
intervals in the class $\cJ$, i.e., with the median set, as it is
defined in the present paper. The purpose of this result is to establish the
equivalence between the definition of multivariate medians in the present work and in
\cite{repold}, for the special case of coordinate-wise partial order, which is considered
there.

\begin{lem}
\label{equalc} Let  $\cI$ and $\cJ$ be  families  of intervals as defined above. Then

\begin{itemize}
\item[(i)] Each interval $I\in \cI$ can be represented as $I= \cap_{J\in \cJ, J\supset I} J$;
\item[(ii)] $\cap_{I\in \cI} I = \cap_{J\in \cJ} J$
\end{itemize}
\end{lem}

\begin{proof}
For $I\in \cI$, let et $S(I)= \bigcap_{J\in \cJ, J\supset I} J$.
Clearly, $I\subset S(I)$. To show that $S(I)\subset I$, take any
$x\not\in I$. Since $I$ is a closed interval, there is another
closed interval $J$ such that
  $I\subset J \subset \{x\}^c$, where both inclusions are strict, and thus $x\not \in S(I)$. This ends the
  proof of (i). To show (ii), note that if $J\in \cJ$, then $J\in \cI$, so $\cJ \subset \cI$. Hence,
\[ \bigcap_{J\in\cJ} J \supset\bigcap_{I\in\cI} I .\]
Conversely, by the part (i), we have that
\[ \bigcap_{I\in \cI} I  = \bigcap_{I} \bigcap_{J\in \cJ, J\supset I} J =\bigcap_{J\in \cJ'}J,\]
where $\cJ'\subset \cJ$, and hence,  we conclude that $
\bigcap_{I\in\cI} I \supset \bigcap_{J\in \cJ} J$.
\end{proof}

{\sc Proof of the Lemma \ref{depf1}.} Evidently, $x\in S_{\alpha}^c$
if and only if $ D(x)<\alpha$, i.e., if and only if there exists a
set $U\in \cU$ such that $x\in U$ and $P(U)<\alpha$. Therefore, if
there are $U\in \cU$ with $P(U)<\alpha$, then
\[ S_{\alpha}^c = \bigcup_{U\in \cU, P(U)<\alpha} U, \quad\mbox{and so,}\quad
S_{\alpha} = \bigcap_{U\in \cU, P(U)<\alpha} U^c ,\] which is
equivalent to the assertion that we wanted to prove. \eop

{\sc Proof of the Theorem \ref{maint}.}  Under $(C_1)$ and if all
sets in  $\cV$ are closed, the set $S_{\alpha}$ is closed for every
$\alpha$, via (\ref{11a}), and hence, the function $D$ is upper
semicontinuous. Under additional conditions $(C_2)$,  we
 will show that there exists at least one $\alpha$ such that
 $S_{\alpha}$ is a nonempty compact set. Indeed, by the assumption, there is $x\in \R^d$ so that
 $D(x)=\alpha_0>0$. On the other hand, by assumption of convergence of $D(x)$ to zero as $\|x\|\tends +\infty$,
 there exists an $R>0$ so that $D(x)<\alpha_0$ for $\|x\|>R$. Therefore, the set $S_{\alpha_0}$ is nonempty and
 norm bounded, and being closed, it is compact. Then all sets $S_{\alpha}$ with $\alpha \geq \alpha_0$ are
 compact,  because $S_{\alpha}\subset S_{\alpha_0}$ for $\alpha \geq \alpha_0$. The intersection of non-empty
 compact nested sets $S_{\alpha}$ is a non-empty compact set, and it is clearly the set on which $D$ reaches its
 maximum. \eop


{\sc Proof of Theorem \ref{equco}. } Suppose that the stated
conditions hold. If $(C_1)$ holds for $\cA$, then $(E')$ implies
that it holds for $\cB$. If $(C_1)$ holds for $\cB$, then it clearly
holds for $\cA$ by $(E'')$.

Let $x\in \R^d$ be fixed. Then by $(E')$, for each $A\in \cA$ that
contains $x$, there exists a $B_A\in \cB$ such that $x\in B_A\subset
A$, and, consequently, $P(A)\geq P(B_A)$. Therefore, \beqn
 D(x; P, \cA) &\geq &  \inf \{ P(B_A)\; |\; x\in B_A\in \cB, A\in \cA\}\\
&\geq & \inf\{ P(B)\; |\; x\in B \in \cB\} = D(x; P, \cB) \eeqn as
the class of all $B_A$ is a subset of the class of all $B\in \cB$
that may contain $x$. On the other hand, by $(E'')$, for each $\eps
>0$ and for each $B\in \cB$ that contains $x$, there exists $A_B\in
\cA$, such that $P(B) \geq P(A_B)-\eps$. Then \beqn
  \inf\{ P(B)\; |\; x\in B \in \cB \} &\geq &
\inf\{ P(A_B)\; |\; x\in A_B, A_B\in \cA, B\in \cB \}-\eps\\
&\geq & \inf\{ P(A)\; |\; x\in A\in \cA\} -\eps = D(x; P, \cA)-\eps,
\eeqn and since $\eps>0$ is arbitrary, we conclude that
\[ D(x; P, \cB) =\inf\{ P(B)\; |\; x\in B \in \cB \} \geq D(x; P, \cA), \]
which ends the proof. \eop

The next Lemma is technical, and we  need it for the proof of
Theorem  \ref{lept}.

\begin{lem}
\label{comlem} Let $P$ be any probability measure on Borel sets of
$\R^d$. Let $K$ be a compact set in $\R^d$ and let $\cA$ be a family
of closed convex subsets of $K$, with $P(A)> \frac{d}{d+1}$ for
every $A\in \cA$. Then the intersection of all sets $A\in \cA$ is a
non-empty compact set.
\end{lem}

\begin{proof}
If $P(A_i) >1-\eps$, $i=1,2,\ldots$, then it is easy to prove by
induction that $P(A_1\cdots A_n) > 1-n\eps$ for $n\geq 2$.
Therefore, under given assumptions, for any $d+1$ sets $A_1,\ldots,
A_n \in \cA$, it holds that $P(A_1\cdots A_{d+1}) > 1-(d+1)\cdot
\frac{1}{d+1}=0$. Hence, every $d+1$ sets of the family $\cA$ have a
non-empty intersection. By Helly's intersection theorem
(\cite[12.12.]{haf}), every finite number of convex sets in $\cA$
have a non-empty intersection. Since $K$ is compact, then all sets
in $\cA$ have a non-empty intersection (see e.g. \cite[Theorem
17.4]{willard}). The intersection is compact since all sets in $\cA$
are compact.
\end{proof}

{\sc Proof of Theorem \ref{lept}.} Let $\delta \in (0,1)$ be fixed.
Assuming that ($C_1$) holds, we will first prove that every compact
convex
set $K\subset \R^d$ with $P(K)=1-\delta>0$ contains a point  $x$
with $D(x; P,\cU)\geq \frac{1-\delta}{d+1}$. Indeed, let $\eps =
\frac{1-\delta}{d+1}$ and suppose, contrary to the statement, that
$D(x; P, \cU) < \eps$ for every $x\in K$, where $K$ is a compact set
with $P(K)=1-\delta >0$. Then (by $(C_1)$), for every $x\in K$ there
exists a $U_x\in \cU$, such that $P(U_x)<\eps$. Clearly, \be
\label{uniu} \bigcup_{x\in K} U_x \supset K. \ee Let $U_x^c=V_x$.
Then $V_x\in \cV$, and by (\ref{uniu}) it follows that \be
\label{intv} \bigcap_{x\in K} (V_x  \cap K) =\emptyset \ee Let us
now define a new probability measure $P^{\ast}$ on $\R^d$, by
\[ P^{\ast} (B) = \frac{P(B\cap K)}{1-\delta},\qquad \mbox{where $B\subset \R^d $ is a Borel set}.\]
For each $x\in K$, we have that $P(V_x) > 1-\eps$, and
\[ P(V_x\cap K) > P(V_x)+P(K)-1> 1-\eps -\delta = \frac{d(1-\delta)}{d+1},\]
hence $P^{\ast} (V_x\cap K) > \frac{d}{d+1}$. Now by Lemma
\ref{comlem}, we conclude that the family of sets $V_x\cap K$ have
non-empty intersection, which contradicts (\ref{intv}). So, the
statement about compact convex sets is proved.

To prove the statement of the Theorem \ref{lept}, note that  the
statement that we already proved yields the condition $(C_2')$, and,
with additional assumption $(C_2'')$, Theorem \ref{maint} is
applicable. By the first part of the proof, each of the sets
\[ S_{n} =\{ x\in \R^d\; |\; D(x; P,\cU)\geq \frac{1-\frac{1}{n}}{d+1}\},\qquad n=1,2,\ldots \]
is non-empty; then  their intersection.
\[ \bigcap_{n=1}^{+\infty} S_n = \{ x\in \R^d\; |\; D(x; P,\cU) \geq \frac{1}{d+1}\}, \]
is also non-empty, by Theorem \ref{maint}. This  ends the proof.
\eop

{\sc Proof of Lemma \ref{condit}.}  We first prove that $(C_2'')$
holds. For a fixed $\eps>0$, and a given probability measure $P$,
 let $B_{1-\eps}$ be a closed ball centered at origin, with $P(B_{1-\eps})>1-\eps$. Then, by assumptions,
 there exists a set $V\in \cV$ such that $B_{1-\eps}\subset V$. By compactness, there exists $r>0$ such
 that all points $x\in V$ satisfy $\|x\|\leq r$. Therefore, all points $x$ with $\| x\| >r$ are in
 $U=V^c$, and, since $P(U) =1-P(V)<\eps$, we conclude that for a given $\eps>0$ there exists $r>0$ so that
 $D(x; P,\cU)< \eps$ for all $x$ with $\|x\|>r$, which proves $(C_2'')$. Then by Theorem \ref{lept},
 the condition $(C_2')$ also holds. \eop

\begin{cor}
\label{cahd} Conditions  $(C_2)$ hold for examples
$2^{\circ}-4^{\circ}$ in \ref{vared}, with any probability measure
$P$.
\end{cor}

\begin{proof}
For examples $2^{\circ}$ and $3^{\circ}$ in \ref{vared} it is
straightforward to check that the assumptions of Lemma \ref{condit}
are satisfied; hence both conditions in  $(C_2)$ hold. For the
halfspace depth in Example $4^{\circ}$, we may apply Theorem
\ref{tequch}, to conclude that the halfspace depth function is the
same as the one based on the family of all compact convex sets,
which is the example $3^{\circ}$.
\end{proof}

{\sc Proof of Theorem \ref{tequch}.} Let $\stackrel{\circ}{H}$ be an
open halfspace from $\stackrel{\circ}{\cH}$, and let $H$ be its
closure. Given any $x\in \stackrel{\circ}{H}$, there exists a closed
halfspace $H_x$ that can be obtained by translation of $H$ in such a
way that the border of $H_x$ contains $x$. Then $H_x\in \cH$ and,
clearly,
\[\stackrel{\circ}{H}= \bigcup_{x\in \stackrel{\circ}{H}}H_x ,\]
which implies condition $(E')$ of Theorem \ref{equco} with
$\cA=\stackrel{\circ}{\cH}$ and $\cB = \cH$. On the other hand, for
any given closed halfspace $H\in \cH$, there exists a sequence of
halfspaces $H_i$, obtained from $H$ by translation, such that
\[\stackrel{\circ}{H_1}\supset \stackrel{\circ}{H_2}\supset \cdots\quad\mbox{and $H=\cap_i \stackrel{\circ}{H_i}
$},\] which is the condition $(E'')$. Therefore, by Theorem
\ref{equco}, \be \label{pf1}
 D(x; P, \cH)=D(x; P, \stackrel{\circ}{\cH}).
\ee Now note that  (\ref{tequch1}) gives condition $(E')$ for $\cA =
\cU$ and $\cB= \stackrel{\circ}{\cH}$ (by taking complements on both
sides); then, as in the proof of Theorem \ref{equco}, we find that
\be \label{pf2}
 D(x; P, \cU) \geq D(x; P, \stackrel{\circ}{\cH}),
 \ee
for every $x\in \R^d$. In the same way, (\ref{tequch2}) gives
condition $(E'')$ for $\cA= \cU$ and $\cB =\cH$, and so, again as in
the proof of \ref{equco}, \be \label{pf3} D(x; P, \cU) \leq D(x; P,
\cH). \ee The statement of the theorem now follows from (\ref{pf1}),
(\ref{pf2}) and (\ref{pf3}). \eop

{\sc Proof of Theorem \ref{jenmedib}.}  Let $Q=Q_{1-\alpha_m}$. Then
for every $\eps>0$, ${\rm Prob}\; (f(X)\leq Q+\eps)
>1-\alpha_m$, and, therefore, the set
\[ V_{\eps} = f^{-1}((-\infty, Q+\eps]) \]
contains the center $C(P,\cU)$. This implies that
\[ f(m)\leq Q+\eps,\quad\mbox{for every $m\in C(P,\cU)$ and every $\eps>0$}.\]
Letting here $\eps\tends 0$, we get (\ref{jenmig}). \eop

{\sc Proof of Theorem \ref{jenmedia}.} From Sections 2 and 3, we
know that, in this case,  the depth function reaches its maximum at
$1/2$; hence, for a given distribution $P$, and the corresponding
random variable $X$, we have that $C(P,\cU)=\{ \Med X\}$, where the
median set is taken with respect to the given partial order
$\preceq$, and $\{ \Med X \} = [\ba_0,\bb_0]$ for some
$\ba_0,\bb_0\in \R^d$.

 Then, by Theorem \ref{jenmedib}, $ f(m) \leq Q_{1/2}$,
for any $m\in \{ \Med X\}$.  If, besides  $Q_{1/2}$, any other
median $M$ of $f(X)$ exists, then we have that $P(f(X)\leq M) =1/2$,
hence the set $V_M = \{ x\; |\; f(x)\leq M\}$ has the probability
$P(V_M)=1/2$. Therefore, $V_M=[\ba,\bb]$ has a non empty
intersection with any interval
$V_{\alpha}=[\ba^{\alpha},\bb^{\alpha}]$ with
$P(V_{\alpha})>\frac{1}{2}$. Then, as in the proof of Lemma
\ref{lemma1}, it follows that $\ba^{\alpha} \preceq \bb$ and $\ba
\preceq \bb^{\alpha}$ for all $\alpha$, which implies, via relations
$\ba_0=\sup_{\alpha} \ba^{\alpha}$ and $\bb_0=\inf_{\alpha}
\bb^{\alpha}$, that
\[ \ba_0 \preceq \bb\quad\mbox{and}\quad \ba \preceq \bb_0,\]
hence, $[\ba ,\bb]\cap [\ba_0,\bb_0]= [\ba_0,\bb_0]\ne \emptyset$. Then the
inequality (\ref{jenmi}) holds with any $m\in [\ba_0,\bb_0]$. \eop

\medskip


\begin{thebibliography}{27}

\bibitem{chakchaud} B. Chakraborty and P. Chaudhuri, {\em On a transformation and re-transformation technique for constructing
  an affine equivariant multivariate median}, Proc. Amer. Math. Soc. {\bf 124 (8)} (1996),
2539--2547.
\bibitem{chentyl} Z. Chen and D. Tyler, {\em The influence function and maximum bias of
{Tukey's} median}, Ann. Statist. {\bf 30 (6)} (2002), 1737--1759.
\bibitem{don1} D.L. Donoho, {\em Breakdown properties of a multivariate location estimators}
PhD thesis, Department of Statistics, Harvard University (1982).
\bibitem{dongas3} D.L. Donoho and M. Gasko, {\em Breakdown properties of location estimates based on halfspace depth
  and projected outlyingness}, Ann. Stat. {\bf 20} (1992), 1803--1827.
\bibitem{hettrand} T.P. Hettmansperger and R.H. Randles, {\em A practical affine equivariant
multivariate median}, Biometrika {\bf 89 (4)} (2002), 851--860.
\bibitem{koshmos} G. Koshevoy and K. Mosler, {\em The {L}orenz zonoid of a
multivariate distribution}, J. Amer. Statist. Assoc. {\bf 91} (1996),
873--882.
\bibitem{liu1} R.Y. Liu, {\em On a notion of simplical depth},
Proc. Nat. Acad. Sci. U.S.A. {\bf 85} (1988), 1732--1734.
\bibitem{liu2} R.Y. Liu, {\em On a notion of data depth based upon random simplices},
Ann. Statist. {\bf 18} (1990), 405--414.
\bibitem{liups1} R.Y. Liu, J.M. Parelius, and K. Singh, {\em Multivariate analysis by data depth:
Descriptive statistics, graphics and inference (with discussion)},
Ann. Stat. {\bf 27} (1999), 783--858.
\bibitem{liusin1} R.Y. Liu and K. Singh, {\em Ordering directional data : concepts of data
depth on circles and spheres}, Ann. Statist. {\bf 20 (3)} (1992),
1468--1484.
\bibitem{masstheo} J.-C. Masse and R. Theorodescu, {\em Halfplane trimming for bivariate
distributions}, J. Multivariate Anal. {\bf 48} (1994), 188--202.
\bibitem{repold} M. Merkle, {\em Jensen's inequality for
medians}, Stat. Prob. Letters {\bf 71} (2005), 277--281.
\bibitem{niinoja} A. Niinimaa and H. Oja, {\em Multivariate median},
Encyclopedia of Statistical Sciences (Update Volume 3). Eds. by
Kotz, S., Johnson, N.L. and Read, C.P., Wiley, New York (1999),
497--505.
\bibitem{nolan} D. Nolan, {\em Asymptotics for multivariate trimming},
Stochastic Processes And Applications {\bf 42} (1992), 157--169.
\bibitem{rocka} Tyrrell R. Rockafellar, {\em Convex Analysis},
Princeton University Press, Princeton (1972).
\bibitem{RouRut} Peter J. Rousseeuw, Ida Ruts, {\em The depth function of a population distribution},
Metrika {\bf 49}(1999), 213--244.
\bibitem{haf} E. Schechter, {\em Handbook of Analysis and Its Foundations},
Academic Press, New York (1997).
\bibitem{small87} C.G. Small, {\em Measures of centrality for multivariate and directional
distributions}, Canadian J. Statistics, {\bf 15} (1987), 31--39.
\bibitem{small4} C.G. Small, {\em A survey of multidimensional medians},
Internat. Stat. Inst. Rev. {\bf 58} (1990), 263--277.
\bibitem{tukey} J.W. Tukey, {\em Mathematics and picturing of data},
Proc. International Congress of Mathematicians, Vancouver,1974 {\bf
volume 2} (1975), 523--531, ISBN 0-919558-04-6.
\bibitem{varzhang} Y. Vardi and C.-H. Zhang, {\em The multivariate {{L}$_{1}$}-median
and associated data depth}, Nat. Acad. Sci. {\bf 97} (2000),
1423--1426.
\bibitem{willard} S. Willard, {\em General Topology},
Addison-Weseley (1970).
\bibitem{zuo1} Y. Zuo, {\em Projection-based depth functions and associated medians},
Ann. Stat. {\bf 31} (2003), 1460--1490.
\bibitem{zuo4} Y. Zuo, {\em Projection based affinne equivariant multivariate
location estimators with the best possible finite sample breakdown
point}, Statist. Sinica {\bf 14: (4)} (2004), 1199--1208.
\bibitem{zuo2} Y. Zuo, {\em Multi-dimensional trimming based on projection depth},
Ann. Stat. {\bf 34 (5)} (2006), 42 pages.
\bibitem{zuocuiy} Y. Zuo, H. Cui, and D. Young, {\em Influence function and maximum bias of
projection depth based estimators}, Ann. Stat. {\bf 32 (1)} (2004),
189--218.
\bibitem{zuoserf1} Y. Zuo and R. Serfling, {\em General notions of statistical depth
function}, Ann. Stat. {\bf 28} (2000a), 461--482.
\bibitem{zuoserf4} Y. Zuo and R. Serfling, {\em Nonparametric notions of multivariate "scatter measure" and "more
  scattered" based on statistical depth functions},
J. Multivariate Anal. {\bf 75} (2000b), 62--78.
\bibitem{zuoserf2} Y. Zuo and R. Serfling, {\em Structural properties and convergence results for contours of sample
  statistical depth functions},
Ann. Stat. {\bf 28} (2000c), 483--499.
\end{thebibliography}

\end{document}